\documentclass[10pt]{iopart}
\usepackage{iopams} 
\usepackage{graphicx}
 
\begin{document}

\title[Higher order fractional derivatives]{Higher order fractional derivatives}

\author{Richard Herrmann}

\address{
GigaHedron, Farnweg 71, D-63225 Langen, Germany
}
\ead{herrmann@gigahedron.com}
\begin{abstract}

Based on the Liouville-Weyl definition of the fractional derivative, a new direct fractional generalization of 
higher order derivatives is presented. It is shown, that the Riesz and Feller derivatives are special cases
of this approach. 

\end{abstract}

%
%
\section{Introduction}
The fractional calculus \cite{f3},\cite{o1} provides a set of axioms and methods to extend the coordinate and corresponding
derivative definitions in a reasonable way from integer order n to arbitrary order $\alpha$:
\begin{equation}
\{ x^n, {\partial^n \over \partial x^n} \} 
\rightarrow
\{ x^\alpha, {\partial^\alpha \over \partial x^\alpha} \}
\end{equation}
The definition of the fractional order derivative is not unique,  several definitions 
e.g. the Riemann, Caputo, Liouville, Weyl, Riesz, Feller, Gr\"unwald  fractional 
derivative definition coexist \cite{caputo}-\cite{pod}. In the last decade, there has been a steadily increasing
interest in applications of the fractional calculus on such different fields of research like mechanics, anomalous 
diffusion or fractional wave equations. 

Most work is dedicated to the special case of the first order derivative operator
($n=1$) replaced by an appropriately chosen fractional derivative operator. This approach is tempting in that sense,
that higher order derivatives my be replaced in a natural way by a consecutive sequence of first order derivatives and 
consequently may be replaced by the corresponding sequence of the fractional extension of the first order derivative. 

Until now, there exists no general approach for a direct fractional extension of higher order derivatives, 
except for ($n=2$). For
that case, Riesz and Feller have derived a fractional generalization of the  second order derivative operator directly. 

In this letter, we will derive a general definition of a direct fractional generalization 
of higher order derivatives. For that purpose, we will first collect the necessary tools, which are currently used
for a fractional generalization of the first- and second order derivative. We will then propose an extension of the 
fractional derivative definition to 
higher order derivatives.

\section{Definitions}     
The Liouville-Weyl fractional integrals of order $0< \alpha <1 $ are defined as
\begin{eqnarray}
I_+^\alpha \phi(x) &=& {1 \over \Gamma(\alpha)} \int_{-\infty}^x (x-\xi)^{\alpha-1}\phi(\xi) d \xi \\
I_-^\alpha \phi(x) &=& {1 \over \Gamma(\alpha)} \int_{x}^{\infty} (\xi-x)^{\alpha-1}\phi(\xi) d \xi 
\end{eqnarray}
The Liouville-Weyl fractional derivatives of order $0< \alpha <1 $ are defined as the left-inverse
operators of the corresponding Liouville-Weyl fractional integrals
\begin{eqnarray}
\label{lw1}
D_+^\alpha \phi(x) &=& I_{+}^{-\alpha}  \phi(x) =   + { \partial \over \partial x}I_+^{1-\alpha} \phi(x) \\
\label{lw2}
D_-^\alpha \phi(x) &=& I_{-}^{-\alpha}  \phi(x) =  -  { \partial \over \partial x}I_-^{1-\alpha} \phi(x) 
\end{eqnarray}
The definitions (\ref{lw1}) and (\ref{lw2}) may be written in an alternative form:
\begin{eqnarray}
\label{Dleft}
D_+^\alpha \phi(x) &=& {\alpha \over \Gamma(1-\alpha)} \int_0^{\infty} {\phi(x)-\phi(x-\xi)\over \xi^{\alpha+1}} d \xi\\
\label{Dright}
D_-^\alpha \phi(x) &=& {\alpha \over \Gamma(1-\alpha)} \int_0^{\infty} {\phi(x)-\phi(x+\xi)\over \xi^{\alpha+1}} d \xi
\end{eqnarray}
which may be derived via
\begin{eqnarray}
D_+^\alpha \phi(x) &=& I_{+}^{-\alpha}  \phi(x) =   + { \partial \over \partial x}I_+^{1-\alpha} \phi(x) \\
 &=& {1 \over \Gamma(1-\alpha)} {\partial \over \partial x} \int_{-\infty}^{x} (x-\xi)^{-\alpha} \phi(\xi) d \xi\\
 &=& {1 \over \Gamma(1-\alpha)} {\partial \over \partial x} \int_0^{\infty} \xi^{-\alpha} \phi(x-\xi) d \xi\\
 &=& {1 \over \Gamma(1-\alpha)}  \int_0^{\infty} \xi^{-\alpha} (-{\partial \over \partial \xi}\phi(x-\xi)) d \xi\\
 &=& {\alpha \over \Gamma(1-\alpha)} \left(  \int_0^{\infty} { \phi(x) \over \xi^{\alpha+1}} d \xi 
                                            -\int_0^{\infty} { \phi(x-\xi) \over \xi^{\alpha+1}} d \xi \right)
\end{eqnarray}
A specific linear combination of the Liouville-Weyl fractional integrals results in the Riesz 
fractional integral $I_R^\alpha$:
\begin{eqnarray}
I_R^\alpha \phi(x) &=&   { I_{+}^{\alpha}+ I_{-}^{\alpha}\over 2 \cos (\alpha \pi/2)}\phi(x) =
\int_{-\infty}^\infty |x-\xi|^{\alpha-1}\phi(\xi) d \xi \quad \alpha>0, \alpha \neq 1,3,5... 
\end{eqnarray}
The Riesz fractional derivative is then given by
\begin{equation}
\label{DR}
D_R^\alpha \phi(x) =   - { D_{+}^{\alpha}+ D_{-}^{\alpha}\over 2 \cos (\alpha \pi/2)}\phi(x)
\end{equation}
or, according to (\ref{Dleft}),(\ref{Dright}):
\begin{equation}
D_R^\alpha \phi(x) =   \Gamma(1+\alpha){\sin(\alpha \pi/2) \over \pi } \int_0^\infty  
   {\phi(x+\xi) - 2 \phi(x) + \phi(x-\xi) \over \xi^{\alpha + 1}} d \xi , \quad 0 < \alpha <2
\end{equation}
Feller has proposed a generalization of the Riesz fractional derivative of the form
\begin{equation}
I_\theta^\alpha \phi(x) =   (c_-(\theta,\alpha)I_+^\alpha + c_+(\theta,\alpha)I_-^\alpha) \phi(x)
\end{equation}
with
\begin{eqnarray}
c_-(\theta,\alpha) &=& {sin((\alpha-\theta)\pi/2) \over \sin(\pi \theta)}\\
c_+(\theta,\alpha) &=& {sin((\alpha+\theta)\pi/2) \over \sin(\pi \theta)}
\end{eqnarray}
The Feller fractional derivative is defined as
\begin{equation}
D_\theta^\alpha \phi(x) =   -\left(c_+(\theta,\alpha)D_+^\alpha + c_-(\theta,\alpha)D_-^\alpha \right) \phi(x)
\end{equation}
Setting $\theta=0$ we obtain
\begin{equation}
c_-(\theta=0,\alpha)=c_+(\theta=0,\alpha) = {1 \over 2 \cos(\alpha \pi/2) }
\end{equation}
which coincides with the definition of the Riesz fractional derivative (\ref{DR}). 

Another special case results for setting $\theta=1$
\begin{equation}
c_-(\theta=1,\alpha)=-c_+(\theta=1,\alpha) = {1 \over 2 \sin(\alpha \pi/2) }
\end{equation}
which leads to the simple form of the fractional derivative:
\begin{eqnarray}
\label{D1}
D_1^\alpha \phi(x) &=& {D_+^\alpha - D_-^\alpha \over 2 \sin(\alpha \pi/2) } \phi(x) \\
                   &=&   \Gamma(1+\alpha) {\cos(\alpha \pi/2) \over \pi } 
\int_0^\infty {\phi(x+\xi) - \phi(x-\xi)  \over \xi^{\alpha+1}} d\xi
\end{eqnarray}
This derivative should be interpreted as the regularized Liouville-Weyl fractional derivative (\ref{DR}).  

Therefore the Feller fractional derivative may be rewritten as a linear combination of $D_1^\alpha$ and $D_0^\alpha$:
\begin{equation}
D_\theta^\alpha \phi(x) =   \left( A_1(\theta,\alpha) (D_+^\alpha - D_-^\alpha)+ A_2(\theta,\alpha) (D_+^\alpha + D_-^\alpha)\right) \phi(x)
\end{equation}
with 
\begin{eqnarray}
A_1(\theta,\alpha) & =&  -\frac{1}{2} \left(c_+(\theta,\alpha) - c_-(\theta,\alpha)\right) = -\frac{1}{2 \sin(\alpha \pi/2)}\sin(\theta \pi/2) \\
A_2(\theta,\alpha) & =&  -\frac{1}{2} \left(c_+(\theta,\alpha) + c_-(\theta,\alpha)\right) = -\frac{1}{2 \cos(\alpha \pi/2)}\cos(\theta \pi/2)
\end{eqnarray}
which finally reads:
\begin{equation}
D_\theta^\alpha \phi(x) = \left(\sin(\theta \pi/2)D_1^\alpha + \cos(\theta \pi/2)D_R^\alpha\right)\phi(x)
\end{equation}
In this form of the Feller fractional derivative the parameter $\theta$ may be interpreted rather as a rotation parameter
instead of a skewness parameter, proposed by other authors. In addition, this form is better suited for a generalization 
to higher order
derivatives, which will be performed in the following section.

\section{The fractional generalization of higher order derivatives}
Using the basic properties of the central differences operators
\begin{eqnarray}
\label{defc1}
\delta_{\frac{1}{2}} \phi(x) & =&  \phi(x + {1 \over 2}\xi)-\phi(x - {1 \over 2}\xi) \\
\label{defc2}
\delta_{1} \phi(x) &=&  {1 \over 2} (\phi(x + \xi)-\phi(x - \xi)) 
\end{eqnarray}
we define the central differences operator $\mathfrak{D}^k $ of order $k$
\begin{equation}
\label{}
\mathfrak{D}^k \phi(x) = \cases {
                               \delta^k_{\frac{1}{2}} \phi(x)  & for k even\\
                               \delta_{1}\delta^{k-1}_{\frac{1}{2}} \phi(x)  & for k odd\\
}
\end{equation}
or explicitely, using (\ref{defc1}) and (\ref{defc2}):
\begin{equation}
\mathfrak{D}^k \phi(x) =    \sum_{n=0}^{2 [(k+1)/2]} a^k_n \phi(x-([(k+1)/2]-n)\xi) 
\end{equation}
with the summation coefficients 
\begin{equation}
a^k_n =  (-1)^n  \cases{\qquad \left(\begin{array}{c} k \\ n \\\end{array} \right) & for k even\\
{1 \over 2}\left[ \left(\begin{array}{c} k-1 \\ n \\\end{array} \right)-\left(\begin{array}{c} k-1 \\ n-2 \\\end{array} \right)\right]  & for k odd\\
}
\end{equation}
The renormalized fractional derivative is then given as:
\begin{equation}
\label{formel}
D^{k;\alpha}
 \phi(x) =  {1 \over N_k} \int_{0}^\infty {d \xi \over \xi^{\alpha+1}}
  \mathfrak{D}^k \phi(x)
\end{equation}
and the normalization factor
\begin{equation}
\fl
N_k = 2 {\Gamma(1-\alpha) \over \alpha} \left( \sum_{n=0}^{[(k+1)/2]} a^k_n (k-n-1)^\alpha \right) 
\cases{
  \cos(\pi \alpha/2) & for k even \\
  \sin(\pi \alpha/2) & for k odd \\
}
\end{equation}
With (\ref{formel}) based on the Liouville definition of the fractional derivative we therefore have given
all fractional derivatives, which extend the ordinary derivative of order $k$:
\begin{equation}
\lim_{\alpha \rightarrow k} {_{k}}D^\alpha = {d^k \over dx^k}
\end{equation}
In addition, for these derivatives the invariance of the scalar product follows:
\begin{equation}
\int_{-\infty}^{\infty}
\left( {_{k}}D^{\alpha *}  f^{*}(x) \right) g(x) dx =
(\pm)^k \int_{-\infty}^{\infty}
f(x)^{*} \left( {_{k}}D^\alpha  g(x) \right)  dx 
\end{equation}
The first four fractional derivative definitions according (\ref{formel}) follow as: 
\begin{eqnarray}
\label{D4}
{_{1}}D^\alpha f(x) &=&  
 \Gamma(1+\alpha) {\cos(\alpha \pi/2) \over \pi } \times \\
&&
\int_0^\infty {f(x+\xi) - f(x-\xi)  \over \xi^{\alpha+1}} d\xi \nonumber\\
&& \qquad \qquad\qquad \qquad\qquad \qquad \qquad  0\leq \alpha < 1 \nonumber \\
{_{2}}D^\alpha f(x) &=&  
 \Gamma(1+\alpha) {\sin(\alpha \pi/2) \over \pi } \times  \\
&&
\int_0^\infty {f(x+\xi) - 2 f(x) + f(x-\xi)  \over \xi^{\alpha+1}} d\xi \nonumber \\
&& \qquad \qquad\qquad \qquad\qquad \qquad \qquad  0\leq \alpha < 2 \nonumber \\
{_{3}}D^\alpha f(x) &=&  
 \Gamma(1+\alpha) {\cos(\alpha \pi/2) \over \pi } {1 \over 2^\alpha -2}\times  \\
&&
\int_0^\infty {-f(x+2 \xi) + 2 f(x+ \xi)-  2 f(x-\xi) + f(x-\xi)  \over \xi^{\alpha+1}} d\xi \nonumber\\
&& \qquad \qquad\qquad \qquad\qquad \qquad \qquad  0\leq \alpha < 3 \nonumber \\
{_{4}}D^\alpha f(x) &=&  
 \Gamma(1+\alpha) {\sin(\alpha \pi/2) \over \pi } {1 \over 2^\alpha -4}\times\\
&&
\fl
\int_0^\infty {-f(x+2 \xi) + 4 f(x+ \xi)- 6 f(x) +4 f(x-\xi) - f(x-2 \xi)  \over  \xi^{\alpha+1}} d\xi \nonumber\\
&& \qquad \qquad\qquad \qquad\qquad \qquad \qquad  0\leq \alpha < 4 \nonumber \\
\end{eqnarray}
These definitions are valid for
$0\leq \alpha < k$. Setting $\alpha>k$
\begin{equation}
{_{k}}D^\alpha = {d^{n k}\over dx^{n k}}{_{k}}D^{\alpha - n k} \qquad ,n \in {\textrm{\bf{N}}}  
\end{equation}
and chosing $n$ so that 
$0 \leq \alpha-n k < k$ 
the definitions given are valid for all $\alpha>0$.

In the same manner the Feller fractional derivative definition may extended to fractional derivatives of higher order.

We introduce hyper sperical coordinates on the unit sphere on
${\bf{\textrm{R}^n}}$:
\begin{eqnarray}
x_1 &=& \cos(\theta_{n-1}) \\
x_2 &=& \sin(\theta_{n-1})\cos(\theta_{n-2})\\
&& ...\nonumber \\
x_{n-1} &=& \sin(\theta_{n-1})\sin(\theta_{n-2})...\cos(\theta_{1})\\
x_{n} &=& \sin(\theta_{n-1})\sin(\theta_{n-2})...\sin(\theta_{1})
\end{eqnarray}
With these coordinates the Feller definition of a fractional derivative may be extended to
\begin{equation}
\label{flrn}
{_\textrm{\tiny{F}}}D_{\{\theta_k\}}^\alpha  = \sum_{k=1}^n x_k \, {_{k}}D^\alpha
\end{equation}
\section{Conclusion}
Based on central differences a generalized fractional derivative of arbitrary order has been propsed.


\end{document}